\documentstyle{article}

\font\boldscriptfont=cmbx7
\font\boldscriptscriptfont=cmbx5

\let\oldbf=\bf

\def\newbf#1{\ifmmode\mathchoice
{\hbox{\oldbf #1}}
{\hbox{\oldbf #1}}
{\hbox{\boldscriptfont #1}}
{\hbox{\boldscriptscriptfont #1}}
\else\oldbf #1\fi}

\reversemarginpar

\reversemarginpar

\def\grayrule{\special{ps:: currentgray 0.9 setgray}
\hbox to 0pt{\hss
\vrule height 0.3 in depth 0.2 in width 0.5 in \hss}
\special{ps:: setgray}}
\def\sqquare#1#2{\square\hbox to 0pt{$#1#2$\hss}}

\newcommand{\blacksquare}{\hbox{\vrule widht 6 pt height 6 pt depth 0pt}}
\newcommand{\Vdash}{\mathrel{{\vrule height 6.9pt depth -0.1pt}\!\vdash}}
\newcommand{\restriction}{|}
\newcommand{\rest}{{\mathord{\restriction}}}

\newcommand{\unif}{{\rm {\bf unif}}\/}

\newcommand{\cf}{{\rm cf}\/}
\newcommand{\dom}{{\rm dom}}

\newcommand{\QED}{\vrule width 6pt height 6pt depth 0pt \vspace{0.1in}}

\newcommand{\forces}{\Vdash}
\newcommand{\F}{{\cal F}}

\newcommand{\mybegin}{\begin{document}}
\newcommand{\myend}{\end{document}}
\newcommand{\myinput}{\input}
\newcommand{\mydocumentstyle}{\documentstyle{book}}
\newcommand{\V}{{\bf V}}
\def\<#1\>{\langle #1 \rangle}
\newcommand{\whi}{\square}
\newcommand{\bla}{\blacksquare}
\newtheorem{theorem}{Theorem}[section]
\newtheorem{lemma}[theorem]{Lemma}
\newtheorem{corollary}[theorem]{Corollary}
\newtheorem{claim}[theorem]{Claim}
\newtheorem{problem}{Problem}
\newtheorem{model}[theorem]{Model}
\newtheorem{definition}[theorem]{Definition}
\newcommand{\lesdot}{\mathrel{\mathord{<}\!\!\raise 
0.8 pt\hbox{$\scriptstyle\circ$}}}
\newcommand{\Proof}{{\sc Proof} \hspace{0.2in}}
\newcommand{\comp}{\circ}

\newcommand{\lless}{\mathord{<}} 
\newcommand{\no}{\|}

\def\lh(#1){|#1|}
\begin{document}

\title{All meager filters may be null}
\author{{\bf Tomek Bartoszynski}\thanks{The author thanks the  Lady Davis
Fellowship 
Trust for full support} \\ 
Boise State University\\
Boise, Idaho\\
and\\
Hebrew University \\
Jerusalem
\and
{\bf Martin Goldstern}\thanks{Supported by the Israel Academy of
Sciences (Basic Research Fund)}\\
Bar Ilan University\\
Ramat Gan, Israel
\and
{\bf Haim Judah}\footnotemark[2]\\
Bar Ilan University\\
Ramat Gan, Israel
\and
{\bf Saharon Shelah}\footnotemark[2]\,\,\thanks % actually, this is not a
					    % thanks but a footnote.
					    % But it seems easier this
					    % way. 
{Publication 434}
\\ 
Hebrew University\\
Jerusalem}
\maketitle
\begin{abstract}
We show that it is consistent with ZFC that all filters which have
the Baire property are Lebesgue measurable. We also show that the existence
of a Sierpinski set implies that there exists a nonmeasurable filter
which has the Baire property.
\end{abstract}
The goal of this paper is to show yet another example of nonduality
between measure and category.

Suppose that ${\cal F}$ is a nonprincipal filter on 
$\omega$. Identify ${\cal F}$ with the set of characteristic functions
of its elements. Under this convention ${\cal F}$ becomes a subset of
$2^{\omega}$ and a question about its topological or
measure-theoretical  properties makes sense.

It has been proved by Sierpinski that every non-principal filter has
either Lebesgue measure zero or is nonmeasurable. Similarly it is
either meager or does not have the Baire property.

In [T] Talagrand  proved that
\begin{theorem}\label{1}
There exists a measurable filter which does not have
the Baire property.$\QED$
\end{theorem}
In fact we have an even stronger result. In [Ba] it is proved that
\begin{theorem}\label{2}
Every measurable filter can be extended to a
measurable filter which does not have the Baire property. $\QED$
\end{theorem}
We show that the dual result is false.
\section{A model where all meager filters are null}
In this section we prove the following theorem:
\begin{theorem}\label{3}
It is consistent with ZFC that every filter which has the Baire property
is measurable.
\end{theorem}
{\sc Proof}
We will use the following more general result:
\begin{theorem}\label{4}
Let $\V \models GCH$ and suppose that $\V[G]$ is a generic extension
extension of $\V$ obtained adding $\omega_{2}$ Cohen reals. 
 Then in
$\V[G]$ for any two sets $A,B \subset 2^{\omega}$ if $A+B=\{a+b:a \in
A, b \in B\}$ is a meager
set then either $A$  or $B$ has measure zero.
\end{theorem}
{\sc Proof} Note that we apply this lemma only for the case $A=B$.
Therefore to simplify the notation we assume that $A=B$. The proof of
the general case is almost the same.
We follow [Bu].

We will use the following notation. Let 
$$Fn(X,2)=\{s : \dom(s) \in [X]^{<\omega} \hbox{ and } \hbox{range}(s)
\subset \{0,1\}\}$$
be the notion of forcing adding $|X|$-many Cohen reals.
For $s \in Fn(X,2)$  let $[s] = \{f \in 2^{X}: s \subset f\}$.

Let $\V \models GCH$ be a model of ZFC and let $G_{\omega_{2}}$ be a
$Fn(\omega_{2},2)$-generic filter over $\V$. Clearly
$c=\bigcup G_{\omega_{2}}$ is a generic sequence of $\omega_{2}$ Cohen
reals and
$\V[c]=\V[G_{\omega_{2}}]$.

Let
$\{F_{n}: n \in \omega \}$ be a sequence of closed, nowhere dense sets
such that $A+A \subseteq \bigcup_{n \in \omega } F_{n}$.
Without loss of  generality we can assume that $\{F_{n}: n \in \omega
\} \in \V$. 

Let $\{a_{\xi}: \xi < \omega_{2}\}$ be an enumeration of all
 elements of $A$.
For every $\xi < \omega_{2}$ let $\dot{a}_{\xi}$ be a name for
$a_{\xi}$. In other words for every $\xi < \omega_{2}$ we have a
countable set $I_{\xi} \subset \omega_{2}$ such that 
$\dot{a}_{\xi}$ is a Borel function from $2^{I_{\xi}}$ into
$2^{\omega}$. Moreover $a_{\xi}$ is the  value of of the function
$\dot{a}_{\xi}$ on Cohen real i.e. $\dot{a}_{\xi}(c\rest
I_{\xi})= a_{\xi}$.
In addition we can find a dense $G_{\delta}$ set $H_{\xi} \subseteq
2^{I_{\xi}}$ such that $\dot{a}_{\xi} \rest H_{\xi}$ is a continuous
function. 

For $\alpha,\xi,\eta < \omega_{2}$ define $\xi
\simeq_{\alpha} \eta$ if
\begin{enumerate}
\item $I_{\xi}$ and $I_{\eta}$ are order isomorphic,
\item the order-isomorphism between $I_{\xi}$ and $I_{\eta}$ transfers
$\dot{a}_{\xi}$ onto $\dot{a}_{\eta}$ and $H_{\xi}$ onto $H_{\eta}$,
\item $I_{\xi} \cap \alpha = I_{\eta} \cap \alpha $.
\end{enumerate}
Notice that for every $\alpha < \omega_{2}$ the relation
$\simeq_{\alpha}$ is an equivalence relation with $\omega_{1}$ many
equivalence classes.
\begin{lemma}
There exists $\alpha^\star < \omega_{2}$ such that
$$\forall \xi,\beta \ \exists \eta \ (\xi \simeq_{\alpha^\star}
\eta \ \& \ I_{\eta} \cap (\beta - \alpha^\star ) = \emptyset) \ .$$
\end{lemma}
{\sc Proof}
For every $\alpha<\omega_2 $ let ${\cal E}_\alpha$ be the set $\{[\xi]_\alpha:
\xi < \omega_2\}$ of $\simeq_\alpha$-equivalence classes.
Let 
$${\cal E}_\alpha^0 = \{ E \in {\cal E}_\alpha: \sup_{\eta \in
E}(\min(I_\eta - \alpha )) < \omega_2\} \hbox{ and}$$
$${\cal E}_\alpha^1 = {\cal E}_\alpha - {\cal E}_\alpha^0 \ .$$
Let 
$$\gamma(\alpha)= \sup_{E \in {\cal E}^0_\alpha}(\sup_{\eta \in
E}(\min(I_\eta -\alpha))) \ .$$
Note that $\gamma(\alpha)<\omega_2$ since $|{\cal
E}_\alpha|<\aleph_1$.

Find $\alpha^\star<\omega_2$ such that $\gamma(\alpha)<\alpha^\star$
for all $\alpha < \alpha^\star$ and $\cf(\alpha^\star)=\omega_1$.
We claim that  $\alpha^\star$ satisfies the statement of the
lemma.

Take any $\xi < \omega_2$ and any $\beta$. If $\beta < \alpha^\star$
or $I_\xi \subseteq \alpha^\star$, then we can choose $\eta = \xi$. So
assume $\beta > \alpha^\star$ and $I_\xi - \alpha^\star \neq
\emptyset$.
There is $\alpha<\alpha^\star$ such that $I_\xi \cap \alpha = I_\xi
\cap \alpha^\star$. Let $E = [\xi]_\alpha$.

{\sc Case 1} $E \in {\cal E}^0_\alpha$. Then
$$\sup_{\eta \in E}(\min(I_\eta - \alpha)) \leq \gamma(\alpha) <
\alpha^\star$$ 
which is a contradiction since $\min(I_\xi - \alpha) \geq
\alpha^\star$ and $\xi \in E$.

{\sc Case 2} $E \not \in {\cal E}^0_\alpha$. So 
$$\sup_{\eta \in E}(\min(I_\eta - \alpha)) = \omega_2$$
hence there is $\eta \in E$ with $\min(I_\eta - \alpha ) \geq \beta$
i.e. $I_\eta \cap (\beta - \alpha)=\emptyset$.

So $I_\xi \cap \alpha^\star =I_\xi \cap \alpha=I_\eta \cap \alpha =
I_\eta \cap \alpha^\star$, where the last equality holds because 
$I_\eta \cap (\alpha^\star - \alpha ) \subseteq I_\eta \cap (\beta -
\alpha)=\emptyset $.
Also $I_\eta \cap (\beta - \alpha^\star) \subseteq I_\eta \cap (\beta
- \alpha)=\emptyset$. $\QED$

Let $\alpha^\star$ be the ordinal from the above lemma. Work in $\V'=\V[c
\rest \alpha^\star ]$.

For every $\xi < \omega_{2}$ define 
$$D_{\xi}=\{s \in Fn(\omega_2-\alpha^\star ,2): \hbox{cl}(\dot{a}_{\xi}([s])) \hbox{
has measure zero } \} \ .$$
\begin{lemma}
$D_{\xi}$ is dense in $Fn(\omega_{2}-\alpha^\star ,2)$ for every $\xi < \omega_{2}$.
\end{lemma}
{\sc Proof}
Notice that it is enough to show that $D_{\xi} \cap Fn(I_\xi -
\alpha^\star,2)$
is dense in
$Fn(I_{\xi}-\alpha^\star,2)$ for $\xi < \omega_{2}$.

Suppose that this fails. Find $\xi<\omega_{2}$ and $s_{0} \in
Fn(I_{\xi}-\alpha^\star ,2 )$ such that
for all $s \supseteq s_0$
 the set $\hbox{cl}(\dot{a}_{\xi}([s]))$
has positive measure.

Using the lemma with $\beta > \sup(I_\xi)$ 
we can find $\eta<\omega_{2}$ such that $\xi
\simeq_{\alpha^\star} \eta$ and $(I_{\xi}-\alpha^\star) \cap
(I_{\eta}-\alpha^\star)=\emptyset $.
Notice that there exists $t_{0} \in Fn(I_{\eta}-\alpha^\star,2)$ (the image
of $s_{0}$ under the isomorphism between $I_{\xi}$ and $I_{\eta}$)
such that for every $t \supseteq t_0$
the set $\hbox{cl}(\dot{a}_{\eta}([t]))$
has positive measure.

Since $s_{0}$ and $t_{0}$ have disjoint domains, $s_{0}\cup t_{0}
\in Fn(\omega_{2}-\alpha^\star,2)$. Find $n \in \omega$ and a condition 
$u \in Fn(\omega_{2}-\alpha^\star,2)$ extending $s_{0}\cup t_{0}$ such
that
$u \forces \dot{a}_{\xi}(\dot{c})+\dot{a}_{\eta}(\dot{c}) \in F_{n}$.
$u$ can be written as  $u_{1}\cup u_{2}\cup u_{3}$ where
$s_{0} \subseteq u_{1} \in Fn(I_{\xi}-\alpha^\star,2)$,
$t_{0} \subseteq u_{2} \in Fn(I_{\eta}-\alpha^\star,2)$ and
$u_{3} \in Fn(\omega_{2}-(I_{\xi} \cup I_{\eta} \cup \alpha^\star),2)$.
By the assumption the sets
$\hbox{cl}(\dot{a}_{\xi}([u_{1}])),
\hbox{cl}(\dot{a}_{\eta}([u_{2}]))$ have positive measure. By
well-known theorem of Steinhaus the set
$\hbox{cl}(\dot{a}_{\xi}([u_{1}]))+
\hbox{cl}(\dot{a}_{\eta}([u_{2}]))$ contains an open set (hence also
$(\hbox{cl}(\dot{a}_{\xi}([u_{1}]))+
\hbox{cl}(\dot{a}_{\eta}([u_{2}])))-F_n$ contains an open set).
Using the fact that $\dot{a}_{\xi}$ and $\dot{a}_{\eta}$ are
continuous functions we can 
find  $u_{1} \subseteq s_{1} \in Fn(I_{\xi}-\alpha^\star,2)$ and
$u_{2} \subseteq t_{1} \in Fn(I_{\eta}-\alpha^\star,2)$ such that 
$(\hbox{cl}(\dot{a}_{\xi}([s_{1}]))+
\hbox{cl}(\dot{a}_{\eta}([t_{1}]))) \cap F_{n} = \emptyset$.
But this is a contradiction since 
$$s_{1}\cup t_{1}\cup u_{3} \forces  
\dot{a}_{\xi}(\dot{c})+\dot{a}_{\eta}(\dot{c}) \not \in F_{n}\ .\QED $$
Notice that for $\xi < \omega_{2}$ 
$$D_{\xi} = \{s \in Fn(I_{\xi}) : \hbox{ there exists  a closed measure
zero set } F \in V'$$
$$ \hbox{ such that } s \forces \dot{a}_{\xi}(\dot{c})
\in F\} \ .$$
Therefore by the above lemma 
$$A \subseteq \bigcup \{F : F \hbox{ is a closed measure zero set
coded in } \V'\} \ .$$
Since $\V$ contains Cohen reals over $\V'$,
the union of all closed measure zero sets coded in $\V'$ has measure
zero in $\V$.
We conclude that $A$ has
measure zero. $\QED$

Let $\F$ be a non-principal filter. Denote by $\F^{c} = \{X \subseteq
\omega: \omega - X \in \F\}$. $\F^{c}$ is an ideal and it is very easy
to see that $\F$ is measurable (has the Baire property) iff $\F^{c}$ is
measurable (has the Baire property).
\begin{lemma}
$\F+\F=\F^{c}$.
\end{lemma}
{\sc Proof}
Suppose that $X,Y \in \F$. Then $\{n : X(n)+Y(n)=0\} \supseteq
X^{-1}(1) \cap Y^{-1}(1) \in \F$.
In general $\F + \cdots + \F$ is equal to $\F$ or $\F^{c}$ depending
whether there is an even or odd number of $\F$'s.

Let $\V \models GCH$ and suppose that $\V[G]$ is a generic 
extension of $\V$ obtained by adding $\omega_{2}$ Cohen reals. 
By the above lemma if $\F$ is a meager filter then $\F^{c} = \F+\F$
is meager. So by \ref{4} $\F$ has measure zero. $\QED$

\section{Filters which are meager and nonmeasurable}
Theorem \ref{3} shows that in order to construct a filter which is
meager and nonmeasurable we need some extra assumptions.

In [T] Talagrand showed that
\begin{theorem}
Suppose that the real line is not the union of $<2^{\aleph_{0}}$ many
measure zero sets. Then there exists a nonmeasurable filter which is
meager.
$\QED$
\end{theorem}
Let $\kappa$ be a regular uncountable cardinal.
Recall that $S$ is a generalized Sierpinski set of size $\kappa$ if  $|S \cap
H< \kappa$ for every null set $H$.
It is clear that  all $S' \subseteq S$ of size $\kappa $ are also
nonmeasurable. 
\begin{theorem}
Assume that there exists a generalized Sierpinski set. Then there exists a
nonmeasurable meager filter.
\end{theorem}
{\sc Proof}
Let $S$ be a generalized Sierpinski set of size $\kappa$. 
Build a sequence $\{x_{\xi}: \xi <
\kappa\} \subset S$ and an elementary chain of  models
$\{M_{\xi}:\xi<\kappa\}$ of size $\kappa$ such that
\begin{enumerate}
\item $\{x_{\xi}:\xi < \alpha \} \subset M_{\alpha}$ for $\alpha <
\kappa$, 
\item $x_{\beta}$ is a random real over $M_{\alpha}$ for
$\beta>\alpha$.
\end{enumerate}
Suppose that $M_{\beta}, x_{\beta}$ are already constructed for $\beta
< \alpha $. Since $S$ is a Sierpinski set $$\bigcup \{S \cap H:
H \hbox{ is a null set coded in } M_{\beta} \hbox{ for } \beta <
\alpha \}$$
has size $<\kappa$. Let $x_{\alpha}$ be any element of $S$
avoiding this set.

Let $X_{\xi}=x_\xi^{-1}(1)$ for $\xi < \kappa $.
Let $\F$ be the filter generated by the family $\{X_{\xi}:\xi <
\kappa \}$. We will show that $\F $ has the required properties.

For $X \subset \omega $ let 
$$d(X) = \lim_{n \rightarrow \infty} \frac{|X \cap n|}{n}$$
if the above limit exists.

By easy induction we show that for $\xi_{1}, \ldots,
\xi_{n}<\kappa $  we have
$d(X_{\xi_{1}} \cap \cdots \cap X_{\xi_{n}}) = 2^{-n}$.
This shows that 
$$\F \subseteq \{X \subset \omega : \liminf_{n \rightarrow \infty} 
\frac{|X \cap n|}{n}>0\}$$
which is
a meager set. To check that $\F$ is nonmeasurable notice that $\F$ contains
the nonmeasurable set $\{x_\xi : \xi < \kappa\}$ .$\QED$

It is an open problem whether one can construct a meager nonmeasurable
filter assuming the existence of a nonmeasurable set of size
$\aleph_{1}$. We only have some partial results.

Let ${\bf b}$ be the size of the smallest unbounded family in
$\omega^{\omega}$ and let $\unif$ be the size of the smallest
nonmeasurable set.

For $X \subseteq \omega $ let $f_{X} \in \omega^{\omega}$ be 
an increasing function enumerating $X$. 
For a filter $\F$ let $\F^{\star} = \{f_{X} : X \in \F\}$.
In [J] it is proved that 
\begin{theorem}\label{haim}
For every filter $\F$,\\
$\F$ has the Baire property iff $\F^{\star}$ is
bounded. $\QED$
\end{theorem}

\begin{theorem}
Suppose that $\unif<{\bf b}$. Then there exists a nonmeasurable filter
which is meager.
\end{theorem}
{\sc Proof}
Let $X \subseteq 2^{\omega}$ be a nonmeasurable set of size $\unif$.
Let $M$ be a model of the same size containing $X$ as a subset.
Then $M \cap 2^\omega $ does not have measure zero, so it is nonmeasurable.
Consider any filter $\F$ such that $M \models \F \hbox{ is an
ultrafilter}$.
$\F$ generates a filter in $\V$ and this filter is meager by 
\ref{haim} and the fact that it is generated by $\unif<{\bf b}$ many
elements. On the other hand $M \models 2^{\omega} = \F \cup
\F^{c}$ and we know that $M \cap 2^{\omega}$ is a nonmeasurable
set. Hence $\F$ is nonmeasurable.$\QED$

The previous theorem depended on the implication:
$$\hbox{If } \F \hbox{ has measure zero then } M \cap 2^\omega
\hbox{ has measure zero.}$$
This implication is not true in general for any set $X \in M$ having
outer measure 1 in $M$ as is showed by the following
example.

{\sc Example} It is consistent with ZFC that there are models $M \subset \V$
such that only {\em some} sets which have outer measure 1 in $M$ have measure
0 in $V$.

Let $\V = {\bf L}[c][\<r_\xi:\xi< \omega_1\>]$ 
where $c$ is a Cohen real over ${\bf L}$ and
$\<r_\xi : \xi < \omega_1\>$ is a sequence of
random reals over ${\bf L}[c]$ (added side by side). 
Let $M={\bf L}[\<r_\xi: \xi < \omega_1\>]$.
Consider the set $X = {\bf L} \cap 2^{\omega}$. It is known that $X$
is a nonmeasurable set in $M$ but $X$ has measure  0 in $\V$.
On the other hand the set $\{r_\xi : \xi < \omega_1\}$ is
nonmeasurable in $\V$.$\QED$

We conclude the paper with a canonical example of a filter which
does not generate an ultrafilter. In other words we have the
following:
\begin{theorem}
Let $M$ be a model for ZFC and let $r$ be a real which does not belong
to $M$. Then there exists a filter $\F$ such that $M \models \F \hbox{
is an ultrafilter}$ but
$$M[r] \models \{X \subseteq \omega : \exists Y \in \F \ Y \subseteq
X\} \hbox{ is not an ultrafilter} \ .$$
\end{theorem}
{\sc Proof}
Let $\{k_{n} : n \in \omega \}$ be a fast increasing sequence of
natural numbers.
Let $T$ be a tree on $2^{<\omega}$ such that:
\begin{enumerate}
\item For $s \in T$ we have  $\lh(s) = k_{n}$ iff $s^{\frown}0 \in T$ and
$s^{\frown}1 \in T$,
\item let $\{s_{1}, \ldots, s_{2^{n}}\}$ be the list of $T \cap
2^{k_{n}}$ in lexicographical order. Then for every $w \subseteq
{\cal P}(2^n)-\{\emptyset, 2^n\}$
there exists $m \in [k_{n}+1,k_{n+1})$ such that
$s_{l}(m)=0$ iff $l \in w$,
\item there is no $m \in \omega $ such that
for all $s \in T \cap 2^{m+1}$ we have $s(m)=0$ or
for all $s \in T \cap 2^{m+1}$ we have $s(m)=1$.
\end{enumerate}

Let $S \subseteq T $ be a subtree of $T$. 
Define 
$$A^{0}_{S} = \{m : \forall s \in S \cap 2^{m+1} \ s(m)=0\} \hbox{
and}$$
$$A^{1}_{S} = \{m : \forall s \in S \cap 2^{m+1} \ s(m)=1\}\ .$$
Let ${\cal J}$ be the ideal generated by sets $\{A^0_S, A^1_S : S
\hbox{ is a perfect subtree of } T\}$.

One can easily verify that all finite subsets of $\omega$ belong to
${\cal J}$. 
\begin{lemma} 
${\cal J}$ is a proper ideal.
\end{lemma}
{\sc Proof}
Let $S_1, \ldots , S_m$ be perfect subtrees of $T$.
Find $n$ sufficiently big so that $|S_j \cap 2^{k_{n}}|>m$ for $j \leq
m$.
Let $s_1, \ldots, s_{2^n}$ be the list of $T \cap 2^{k_{n}}$ in
lexicographical ordering. Let $w_{1}, \ldots, w_{m}$ be such that 
$S_j \cap 2^{k_n}= \{s_i : i \in w_j\}$ for $j \leq m$.
Let $w = \{\min(w_1), \ldots, \min(w_m)\}$. 
Then for all $j, \ w_j \not \subseteq w$ and $w_j \cap w \neq \emptyset$.
By the definition of $T$ there is $k<k_n$ such that $w = \{l :
s_l(k)=0\}$.
By the property of $w$ for every $j \leq m$ there exist $s^0,s^1 \in
S_j \cap 2^{k_{n}}$ such that $s^0(k)=0$ and $s^1(k)=1$.
Therefore $k \not \in A^0_{S_1} \cup A^1_{S_1} \cup \cdots \cup A^0_{S_m}
\cup A^1_{S_m}$.

Let $\F$ be any ultrafilter in $M$ extending the filter $\{\omega - X
: X \in {\cal J}\}$.
 Let $r$ be a real which does not belong to $M$. Without loss of
generality we can assume that $r$ is a branch through $T$.

Assume that $\F$ generates an ultrafilter and let $X_{r} = \{n :
r(n)=1\}$.
We can  assume that there exists an element $X \in \F$ such that 
$X \subseteq X_r$.
Let $S=\{s \in T : \forall k \in X \ (\lh(s) > k \ \rightarrow
s(k)=1)\}$.
Clearly $r$ is a branch through $S$. But in that case
$S$ contains a perfect subtree $S_1 \subseteq S$ (since it contains a
new branch).
Therefore $X \subseteq A^1_{S_1} \in {\cal J}$. Contradiction. $\QED$


\begin{thebibliography}{tomek}
\bibitem[Ba]{Ba} T. Bartoszynski {\em On the structure of the filters
on a countable set} to appear
\bibitem[Bu]{Bu} M. Burke {\em notes of June 17, 1989}
\bibitem[J]{J} H. Judah {\em Unbounded filters on $\omega$}, in {\bf
Logic Colloquium 1987}
\bibitem[T]{T} M. Talagrand {\em Compacts de
fonctions mesurables et filtres nonmesurables},
{\bf Studia Mathematica}, T.LXVII, 1980.
\end{thebibliography}
\end{document}